\documentclass[12pt,reqno]{amsart}
\usepackage{amssymb}
\usepackage{amsxtra}
\usepackage[mathscr]{eucal}

\newcommand{\eop}{\hfill$\square$}

\theoremstyle{plain}
\newtheorem{Thm}{Theorem}

\newtheorem{Lem}{Lemma}

\theoremstyle{definition}

\theoremstyle{remark}

\begin{document}

\title[David M.~Bradley]{A $q$-Analog of Euler's Decomposition Formula
for the Double Zeta Function}

\date{\today}

\author{David~M. Bradley}
\address{Department of Mathematics \& Statistics\\
         University of Maine\\
         5752 Neville Hall
         Orono, Maine 04469-5752\\
         U.S.A.}
\email[]{bradley@math.umaine.edu, dbradley@member.ams.org}



\begin{abstract}
The double zeta function was first studied by Euler in response to
a letter from Goldbach in 1742.  One of Euler's results for this
function is a decomposition formula, which expresses the product
of two values of the Riemann zeta function as a finite sum of
double zeta values involving binomial coefficients.  In this note,
we establish a $q$-analog of Euler's decomposition formula. More
specifically, we show that Euler's decomposition formula can be
extended to what might be referred to as a ``double $q$-zeta
function'' in such a way that Euler's formula is recovered in the
limit as $q$ tends to 1.
\end{abstract}

\maketitle

\interdisplaylinepenalty=500

\section{Introduction}\label{sect:Intro}
The Riemann zeta function is defined for $\Re(s)>1$ by
\begin{equation}\label{Riemann}
   \zeta(s) := \sum_{n=1}^\infty \frac{1}{n^s}.
\end{equation}
Accordingly,
\begin{equation}\label{doublezeta}
   \zeta(s,t)
   := \sum_{n=1}^\infty \frac{1}{n^s}\sum_{k=1}^{n-1}\frac{1}{k^t},
      \quad \Re(s)>1, \quad \Re(s+t)>2,
\end{equation}
is known as the double zeta function.  The
sums~\eqref{doublezeta}, and more generally those of the form
\begin{equation}\label{MzvDef}
   \zeta(s_1,s_2,\dots,s_m) := \sum_{k_1>k_2>\cdots>k_m>0}\;
   \prod_{j=1}^m \frac{1}{k_j^{s_j}},
   \quad \sum_{j=1}^n \Re(s_j) > n,\quad n=1,2,\dots,m,
\end{equation}
have attracted increasing attention in recent years; see
eg.~\cite{BBB,BBBLa,BBBLc,BowBrad1,BowBrad3,BowBradRyoo,Prtn,DBqKarl,BK1,LeM}.
The survey articles~\cite{BowBradSurvey,Cartier,Wald,Wald2,Zud}
provide an extensive list of references. In~\eqref{MzvDef} the sum
is over all positive integers $k_1,\dots,k_m$ satisfying the
indicated inequalities.  Note that with positive integer
arguments, $s_1>1$ is necessary and sufficient for convergence.

The problem of evaluating sums of the form~\eqref{doublezeta} for
integers $s>1$, $t>0$ seems to have been first proposed in a
letter from Goldbach to Euler~\cite{LE2} in 1742. (See
also~\cite{LE,Goldbach} and~\cite[p.\ 253]{Berndt1}.)  Among other
results for~\eqref{doublezeta}, Euler proved that if $s-1$ and
$t-1$ are positive integers, then the decomposition formula
\begin{equation}\label{Decomp}
   \zeta(s)\zeta(t)
   = \sum_{a=0}^{s-1}\binom{a+t-1}{t-1}\zeta(t+a,s-a)
   + \sum_{a=0}^{t-1}\binom{a+s-1}{s-1}\zeta(s+a,t-a)
\end{equation}
holds.   A combinatorial proof of Euler's decomposition
formula~\eqref{Decomp} based on the Drinfel'd integral
representations~\cite{BBB,BBBLa,BBBLc,BowBradSurvey,BowBrad3}
\begin{gather}\label{iterints}
\begin{split}
   \zeta(s) &= \int\limits_{1>x_1>\cdots>x_s>0}
               \bigg(\prod_{i=1}^{s-1} \frac{dx_i}{x_i}\bigg)
               \frac{dx_s}{1-x_s},\\
   \zeta(s,t) &= \int\limits_{1>x_1>\cdots>x_{s+t}>0}
          \bigg(\prod_{i=1}^{s-1} \frac{dx_i}{x_i}\bigg)
          \frac{dx_s}{1-x_s}
          \bigg(\prod_{i=s+1}^{s+t-1} \frac{dx_i}{x_i}\bigg)
          \frac{dx_{s+t}}{1-x_{s+t}},
\end{split}
\end{gather}
and the shuffle multiplication rule satisfied by such integrals is
given in~\cite[eq.\ (10)]{BBBLc}.  It is of course well-known
that~\eqref{Decomp} can also be proved algebraically by summing
the partial fraction decomposition~\cite[p.\
48]{Niels}~\cite[Lemma 3.1]{Mark}
\begin{equation}\label{parfrac}
   \frac{1}{x^s(c-x)^t}
   = \sum_{a=0}^{s-1}\binom{a+t-1}{t-1}\frac{1}{x^{s-a}c^{t+a}}
   +
   \sum_{a=0}^{t-1}\binom{a+s-1}{s-1}\frac{1}{c^{s+a}(c-x)^{t-a}}
\end{equation}
over appropriately chosen integers $x$ and $c$.  (See
eg.~\cite{BBG}.)

A $q$-analog of~\eqref{MzvDef} was independently introduced
in~\cite{DBqMzv,OkudaYoshihiro,Zhao} as
\begin{equation}
   \zeta[s_1,s_2,\dots,s_m] := \sum_{k_1>k_2>\cdots >k_m>0}
   \; \prod_{j=1}^m \frac{q^{(s_j-1)k_j}}{[k_j]_q^{s_j}},
   \label{qMzvDef}
\end{equation}
where
\[
   [k]_q := \sum_{j=0}^{k-1} q^j = \frac{1-q^k}{1-q},
   \qquad 0<q<1.
\]
Observe that we now have
\[
   \zeta(s_1,\dots,s_m) = \lim_{q\to 1-} \zeta[s_1,\dots,s_m],
\]
so that~\eqref{qMzvDef} represents a generalization
of~\eqref{MzvDef}.  In this note, we establish a $q$-analog of
Euler's decomposition formula~\eqref{Decomp}.

\section{Main Result}\label{sect:main}

Our $q$-analog of Euler's decomposition formula naturally requires
only the $m=1$ and $m=2$ cases of~\eqref{qMzvDef}; specifically
the $q$-analogs of~\eqref{Riemann} and~\eqref{doublezeta} given by
\begin{equation}\label{qzetas}
    \zeta[s] = \sum_{n>0} \frac{q^{(s-1)n}}{[n]_q^s}
    \quad \text{and} \quad
    \zeta[s,t] =
  \sum_{n>k>0}\frac{q^{(s-1)n}q^{(k-1)t}}{[n]_q^s[k]_q^t}.
\end{equation}
We also define, for convenience, the sum
\begin{equation}\label{phi}
   \varphi[s] := \sum_{n=1}^\infty \frac{(n-1)q^{(s-1)n}}{[n]_q^s}
   = \sum_{n=1}^\infty \frac{nq^{(s-1)n}}{[n]_q^s} - \zeta[s].
\end{equation}

We can now state our main result.
\begin{Thm}\label{thm:qDecomp} If $s-1$ and $t-1$ are positive
integers, then
\begin{align*}
   \zeta[s]\zeta[t]
   &= \sum_{a=0}^{s-1}\,\sum_{b=0}^{s-1-a}
     \binom{a+t-1}{t-1}\binom{t-1}{b}(1-q)^b\,\zeta[t+a,s-a-b]\\
   &+ \sum_{a=0}^{t-1}\,\sum_{b=0}^{t-1-a}
     \binom{a+s-1}{s-1}\binom{s-1}{b}(1-q)^b\,\zeta[s+a,t-a-b]\\
   &- \sum_{j=1}^{\min(s,t)}\frac{(s+t-j-1)!}{(s-j)!\,(t-j)!}\cdot
      \frac{(1-q)^j}{(j-1)!}\,\varphi[s+t-j].
\end{align*}
\end{Thm}
Observe that the limiting case $q=1$ of Theorem~\ref{thm:qDecomp}
reduces to Euler's decomposition formula~\eqref{Decomp}.

\section{A Differential Identity}\label{sect:diff}
Our proof of Theorem~\ref{thm:qDecomp} relies on the following
identity.

\begin{Lem}\label{lem:diff} Let $s$ and $t$ be positive
integers, and let $x$ and $y$ be non-zero real numbers.  Then for
all real $q$,
\begin{align*}
   \frac{1}{x^s y^t}
   &= \sum_{a=0}^{s-1}\,\sum_{b=0}^{s-1-a}
      \binom{a+t-1}{t-1}\binom{t-1}{b}
      \frac{(1-q)^b(1+(q-1)y)^a(1+(q-1)x)^{t-1-b}}
           {x^{s-a-b}(x+y+(q-1)xy)^{t+a}}\\
   &+ \sum_{a=0}^{t-1}\,\sum_{b=0}^{t-1-a}
      \binom{a+s-1}{s-1}\binom{s-1}{b}
      \frac{(1-q)^b(1+(q-1)x)^a(1+(q-1)y)^{s-1-b}}
           {y^{t-a-b}(x+y+(q-1)xy)^{s+a}}\\
   &- \sum_{j=1}^{\min(s,t)}\frac{(s+t-j-1)!}{(s-j)!\,(t-j)!}
      \cdot\frac{(1-q)^j}{(j-1)!}\cdot
      \frac{(1+(q-1)y)^{s-j}(1+(q-1)x)^{t-j}}{(x+y+(q-1)xy)^{s+t-j}}.
\end{align*}
\end{Lem}

\noindent{\bf Proof.} Apply the partial differential operator
\[
   \frac{1}{(s-1)!}\bigg(-\frac{\partial}{\partial x}\bigg)^{s-1}
   \frac{1}{(t-1)!}\bigg(-\frac{\partial}{\partial y}\bigg)^{t-1}
\]
to both sides of the identity
\[
   \frac{1}{xy}
 = \frac{1}{x+y+(q-1)xy}\bigg(\frac{1}{x}+\frac{1}{y}+q-1\bigg).
\]
\eop

Observe that when $q=1$, Lemma~\ref{lem:diff} reduces to the
identity
\[
   \frac{1}{x^sy^t}
 = \sum_{a=0}^{s-1}\binom{a+t-1}{t-1}\frac{1}{x^{s-a}(x+y)^{t+a}}
 + \sum_{a=0}^{t-1}\binom{a+s-1}{s-1}\frac{1}{(x+y)^{s+a}y^{t-a}},
\]
from which the partial fraction identity~\eqref{parfrac} (proved
by induction in~\cite{Mark}) trivially follows.

\section{Proof of Theorem~\ref{thm:qDecomp}}

First, observe that if $s>1$ and $t>1$, then from~\eqref{qzetas},
\[
   \zeta[s]\zeta[t]
  = \sum_{n=1}^\infty\, \sum_{u+v=n}
    \frac{q^{(s-1)u}}{[u]_q^s} \cdot \frac{q^{(t-1)v}}{[v]_q^t},
\]
where the inner sum is over all positive integers $u$ and $v$ such
that $u+v=n$.  Next, apply Lemma~\ref{lem:diff} with $x=[u]_q$,
$y=[v]_q$, noting that then
\[
  1+(q-1)x = q^u, \qquad 1+(q-1)y = q^v, \qquad
  x+y+(q-1)xy = [u+v]_q.
\]
After interchanging the order of summation, there comes
\begin{align*}
  \zeta[s]\zeta[t]
  &= \sum_{a=0}^{s-1}\,\sum_{b=0}^{s-1-a}
     \binom{a+t-1}{t-1}\binom{t-1}{b}(1-q)^b S[s,t,a,b]\\
  &+ \sum_{a=0}^{t-1}\,\sum_{b=0}^{t-1-a}
     \binom{a+s-1}{s-1}\binom{s-1}{b}(1-q)^b S[t,s,a,b]\\
  &- \sum_{j=1}^{\min(s,t)}\frac{(s+t-j-1)!}{(s-j)!\,(t-j)!}\cdot
     \frac{(1-q)^j}{(j-1)!}\, T[s,t,j],
\end{align*}
where
\begin{align*}
   S[s,t,a,b]
  &= \sum_{n=1}^\infty\,\sum_{u+v=n}
     \frac{q^{(s-1)u} q^{(t-1)v} q^{(t-1-b)u} q^{av}}
          {[u]_q^{s-a-b} [u+v]_q^{t+a}}
  = \sum_{n=1}^\infty\, \sum_{u+v=n}
     \frac{q^{(t+a-1)(u+v)} q^{(s-a-b-1)u}}
          {[u+v]_q^{t+a} [u]_q^{s-a-b}}\\
  &= \sum_{n=1}^\infty \frac{q^{(t+a-1)n}}{[n]_q^{t+a}}
     \sum_{u=1}^{n-1} \frac{q^{(s-a-b-1)u}}{[u]_q^{s-a-b}}\\
  &= \zeta[t+a,s-a-b]
\end{align*}
and
\begin{align*}
   T[s,t,j]
  &= \sum_{n=1}^\infty \,\sum_{u+v=n}
     \frac{q^{(s-1)u} q^{(t-1)v} q^{(t-j)u} q^{(s-j)v}}
          {[u+v]_q^{s+t-j}}
   = \sum_{n=1}^\infty\, \sum_{u+v=n}
     \frac{q^{(s+t-j-1)(u+v)}}{[u+v]_q^{s+t-j}}\\
  &= \varphi[s+t-j].
\end{align*}
\eop

\section{Final Remarks}\label{sect:final}
In~\cite{Zhao}, Zhao also gives a formula for the product
$\zeta[s]\zeta[t]$.  However, Zhao's formula is considerably more
complicated than ours, as it is derived based on the $q$-shuffle
rule~\cite{BowBradSurvey,DBqMzv} satisfied by the Jackson
$q$-integral analogs of the representations~\eqref{iterints}.  Of
course, we also have the very simple $q$-stuffle~\cite{DBqMzv}
formula
$\zeta[s]\zeta[t]=\zeta[s,t]+\zeta[t,s]+\zeta[s+t]+(1-q)\zeta[s+t-1]$.

\end{document}